\documentclass[12pt]{article}
\usepackage{amsfonts}
\usepackage{amsmath}
\usepackage{graphicx}

\input tcilatex
\input tcilatex
\begin{document}

\author{Ozlem Ersoy \thanks{%
Ozlem Ersoy, Telephone: +90 222 2393750, Fax: +90 222 2393578, E-mail:
ozersoy@ogu.edu.tr} and Idris Dag \\
%EndAName
Mathematics-Computer Department,\\
Eski\c{s}ehir Osmangazi University, 26480, Eski\c{s}ehir, Turkey.$^{1}$}
\title{An Exponential Cubic B-spline Finite Element Method for Solving the
Nonlinear Coupled Burger Equation}
\maketitle

\begin{abstract}
The exponential cubic B-spline functions together with Crank Nicolson are
used to solve numerically the nonlinear coupled Burgers' equation using
collocation method. This method has been tested by three different problems.
The proposed scheme is compared with some existing methods. We have noticed
that proposed scheme produced a highly accurate results.
\end{abstract}

\section{Introduction}

\qquad The purpose of this paper is to apply the exponential B-spline
collocation method to the coupled Burgers equation system. The Coupled
Burger equation system in the following form

\begin{equation}
\begin{array}{c}
\dfrac{\partial U}{\partial t}-\dfrac{\partial ^{2}U}{\partial x^{2}}+k_{1}U%
\dfrac{\partial U}{\partial x}+k_{2}(UV)_{x}=0 \\ 
\dfrac{\partial V}{\partial t}-\dfrac{\partial ^{2}V}{\partial x^{2}}+k_{1}V%
\dfrac{\partial V}{\partial x}+k_{3}(UV)_{x}=0%
\end{array}
\label{CB1}
\end{equation}%
where $k_{1}$, $k_{2}$ and $k_{3}$ are reel constants and subscripts $x$ and 
$t$ denote differentiation, $x$ distance and $t$ time, is considered.
Boundary conditions%
\begin{equation}
\text{\ }%
\begin{array}{l}
U(a,t)=f_{1}(a,t),\text{ }U(b,t)=f_{2}(b,t) \\ 
V(a,t)=g_{1}(a,t),\text{ }V(b,t)=g_{2}(b,t),\text{ }t>0%
\end{array}
\label{CB2}
\end{equation}%
and initial conditions%
\begin{eqnarray}
U(x,0) &=&f(x)  \label{CB3} \\
V(x,0) &=&g(x),\text{ }a\leq x\leq b  \notag
\end{eqnarray}%
will be decided in the later sections according to test problem.

Various methods are used solve the nonlinear coupled Burgers' equation
numerically; which is suggested by \cite{1995} firstly. Fourth order
accurate compact ADI scheme \cite{1999}, A chebyshev spectral collocation
method \cite{2008}, A meshfree technique \cite{2009}, the Fourier
pseudospectral method \cite{2010}, the generalized two-dimensional
differential transform method \cite{2011}, cubic B-spline collocation method 
\cite{2011iki}, generalized differential quadrature method \cite{2011uc}, a
robust technique for solving optimal control of coupled Burgers' equations 
\cite{2011dort}, a differential quadrature method \cite{2012}, Galerkin
quadratic B-spline finite element method \cite{2012iki}, a fully implicit
finite-difference method \cite{2013}, a composite numerical scheme based on
finite difference \cite{2014}, an implicit logarithmic finite difference
method \cite{2014iki}, modified cubic B-spline collocation method \cite%
{2014uc} was applied to obtain numerical solution of nonlinear Coupled
Burgers system. There are not many articles about exponential cubic B-spline
method for the solving nonlinear differantial equation system.

The exponential splines and exponential B-splines are defined as a
generalization of the well-known splines and B-splines by McCartin\cite%
{mc1,mc2,mc3}. He has also showed a reliable algorithm by using the
exponential spline functions to solve the hyperbolic conservation laws
McCartin\cite{mc4}. The exponential B-splines include a free parameter which
cause to have different bell like piece wise polynomial. The best free
parameter is determined for the exponential B-spline functions for solving
the differential equations. The use of the exponential B-splines is not as
common as the well known B-splines. There are a few studies existing to use
exponential B-splines for build up numerical methods. The singularly
perturbed boundary value problem has been solved based on the collocation
methods with the exponential B-splines \cite{s1,s2,s3}. Very recently,
Exponential B-spline collocation method is applied to compute numerical
solution of the convection diffusion equation\cite{cv}.

The paper is organized as follows. In Section 2, some details about
exponential cubic B-spline collocation method are provided. In Section 3,
the initial states are documented. In section 4 , numerical results for
three different problems and some related figures are given in order to show
the efficiency as well as the accuracy of the proposed method. Finally,
conclusions are followed in Section 5.

\section{Exponential Cubic B-spline Collocation Method}

\qquad Let $\pi $ be partition of the problem domain $[a,b]$ defined at the
knots 
\begin{equation*}
\pi :a=x_{0}<x_{1}<\ldots <x_{N}=b
\end{equation*}%
with mesh spacing $h=(b-a)/N.$ The exponential B-splines, $B_{i}${}$(x),$
with knots at the points of $\pi $ can be defined as

\begin{equation}
B_{i}(x)=\left \{ 
\begin{array}{ll}
b_{2}\left( \left( x_{i-2}-x\right) -\dfrac{1}{p}\left( \sinh (p\left(
x_{i-2}-x\right) )\right) \right)  & \left[ x_{i-2},x_{i-1}\right] , \\ 
a_{1}+b_{1}\left( x_{i}-x\right) +c_{1}\exp \left( p\left( x_{i}-x\right)
\right) +d_{1}\exp \left( -p\left( x_{i}-x\right) \right)  & \left[
x_{i-1},x_{i}\right] , \\ 
a_{1}+b_{1}(x-x_{i})+c_{1}\exp \left( p\left( x-x_{i}\right) \right)
+d_{1}\exp \left( -p\left( x-x_{i}\right) \right)  & \left[ x_{i},x_{i+1}%
\right] , \\ 
b_{2}\left( (x-x_{i+2})-\dfrac{1}{p}(\sinh \left( p\left( x-x_{i+2}\right)
\right) )\right)  & \left[ x_{i+1},x_{i+2}\right] , \\ 
0 & \text{otherwise.}%
\end{array}%
\right.   \label{CB4}
\end{equation}%
where%
\begin{equation*}
\begin{array}{l}
a_{1}=\dfrac{phc}{phc-s},\text{ }b_{1}=\dfrac{p}{2}\left( \dfrac{c(c-1)+s^{2}%
}{(phc-s)(1-c)}\right) ,\text{ }b_{2}=\dfrac{p}{2(phc-s)}, \\ 
\\ 
c_{1}=\dfrac{1}{4}\left( \dfrac{\exp (-ph)(1-c)+s(\exp (-ph)-1)}{(phc-s)(1-c)%
}\right) , \\ 
\\ 
d_{1}=\dfrac{1}{4}\left( \dfrac{\exp (ph)(c-1)+s(\exp (ph)-1)}{(phc-s)(1-c)}%
\right) .%
\end{array}%
\end{equation*}%
and $s=\sinh (ph),$ $c=\cosh (ph),$ $p$ is a free parameter. When $p=1$,
graph of the exponential cubic B-splines over the interval $[0.1]$ is
depicted in Fig. 1.%
\begin{equation*}
\begin{array}{c}
\FRAME{itbpF}{2.4431in}{2.4431in}{0in}{}{}{fig1.bmp}{\special{language
"Scientific Word";type "GRAPHIC";maintain-aspect-ratio TRUE;display
"USEDEF";valid_file "F";width 2.4431in;height 2.4431in;depth
0in;original-width 3.2197in;original-height 3.2197in;cropleft "0";croptop
"1";cropright "1";cropbottom "0";filename 'Fig1.bmp';file-properties
"XNPEU";}} \\ 
\text{Fig.1: Exponential cubic B-splines for }p=1\text{ over the interval }%
[0.1]%
\end{array}%
\end{equation*}

An additional knots outside the problem domain, positioned at $x_{-1}<x_{0}$
and $x_{N}<x_{N+1}$ are necessary to define all exponential splines. So that 
$\{B_{-1}(x),B_{0}(x),\cdots ,B_{N+1}(x)\}$ forms a basis for the functions
defined over the interval. Each sub interval $[x_{i},x_{i+1}]$ is covered by
four consecutive exponential B-splines. The exponential B-splines and its
first and second derivatives vanish outside its support interval $%
[x_{i-2},x_{i+2}].$ Each basis function $B_{i}(x)$ is twice continuously
differentiable. The values of $B_{i}(x),$ $B_{i}^{^{\prime }}(x)$ and $%
B_{i}^{^{\prime \prime }}(x)$ at the knots $x_{i}$ can be computed from Eq.(%
\ref{CB4}) are shown Table 1.

\begin{equation*}
\begin{tabular}{l}
Table 1: Values of $B_{i}(x)$ and its principle two \\ 
derivatives at the knot points \\ 
$%
\begin{tabular}{|l|ccccc|}
\hline
$x$ & $x_{i-2}$ & $x_{i-1}$ & $x_{i}$ & $x_{i+1}$ & $x_{i+2}$ \\ \hline
$B_{i}$ & $0$ & $\dfrac{s-ph}{2(phc-s)}$ & $1$ & $\dfrac{s-ph}{2(phc-s)}$ & $%
0$ \\ 
$B_{i}^{^{\prime }}$ & $0$ & $\dfrac{p(1-c)}{2(phc-s)}$ & $0$ & $\dfrac{%
p(c-1)}{2(phc-s)}$ & $0$ \\ 
$B_{i}^{^{\prime \prime }}$ & $0$ & $\dfrac{p^{2}s}{2(phc-s)}$ & $-\dfrac{%
p^{2}s}{phc-s}$ & $\dfrac{p^{2}s}{2(phc-s)}$ & $0$ \\ \hline
\end{tabular}%
$%
\end{tabular}%
\end{equation*}

An approximate solution $U_{N}(x,t)$ and $V_{N}(x,t)$ to the analytical
solution $U(x,t)$ and $V(x,t)$ can be assumed of the forms

\begin{equation}
U_{N}(x,t)=\sum_{i=-1}^{N+1}\delta _{i}B_{i}(x),\text{ }V_{N}(x,t)=%
\sum_{i=-1}^{N+1}\phi _{i}B_{i}(x)  \label{CBu1}
\end{equation}%
where $\delta _{i}$ are time dependent parameters to be determined from the
collocation method. The first and second derivatives also can be defined by 
\begin{eqnarray}
U_{N}^{\prime }(x,t) &=&\sum_{i=-1}^{N+1}\delta _{i}B_{i}^{\prime }(x),\text{
}V_{N}^{\prime }(x,t)=\sum_{i=-1}^{N+1}\phi _{i}B_{i}(x)  \label{CBu2} \\
U_{N}^{\prime \prime }(x,t) &=&\sum_{i=-1}^{N+1}\delta _{i}B_{i}^{\prime
\prime }(x),\text{ }V_{N}^{\prime \prime }(x,t)=\sum_{i=-1}^{N+1}\phi
_{i}B_{i}(x)  \notag
\end{eqnarray}

Using the Eq. (\ref{CBu1}), (\ref{CBu2}) and Table 1, we see that the nodal
values $U_{i},$ $V_{i}$, their first derivatives $U_{i}^{\prime },$ $%
V_{i}^{\prime }$ and second derivatives $U_{i}^{\prime \prime },$ $%
V_{i}^{\prime \prime }$ at the knots are given in terms of parameters by the
following relations%
\begin{equation}
\begin{array}{c}
\begin{tabular}{l}
$U_{i}=U(x_{i},t)=\dfrac{s-ph}{2(phc-s)}\delta _{i-1}+\delta _{i}+\dfrac{s-ph%
}{2(phc-s)}\delta _{i+1},$ \\ 
$U_{i}^{\prime }=U^{\prime }(x_{i},t)=\dfrac{p(1-c)}{2(phc-s)}\delta _{i-1}+%
\dfrac{p(c-1)}{2(phc-s)}\delta _{i+1}$ \\ 
$U_{i}^{\prime \prime }=U^{\prime \prime }(x_{i},t)=\dfrac{p^{2}s}{2(phc-s)}%
\delta _{i-1}-\dfrac{p^{2}s}{phc-s}\delta _{i}+\dfrac{p^{2}s}{2(phc-s)}%
\delta _{i+1}.$%
\end{tabular}
\\ 
\begin{tabular}{l}
$V_{i}=V(x_{i},t)=\dfrac{s-ph}{2(phc-s)}\phi _{i-1}+\phi _{i}+\dfrac{s-ph}{%
2(phc-s)}\phi _{i+1},$ \\ 
$V_{i}^{\prime }=V^{\prime }(x_{i},t)=\dfrac{p(1-c)}{2(phc-s)}\phi _{i-1}+%
\dfrac{p(c-1)}{2(phc-s)}\phi _{i+1}$ \\ 
$V_{i}^{\prime \prime }=V^{\prime \prime }(x_{i},t)=\dfrac{p^{2}s}{2(phc-s)}%
\phi _{i-1}-\dfrac{p^{2}s}{phc-s}\phi _{i}+\dfrac{p^{2}s}{2(phc-s)}\phi
_{i+1}.$%
\end{tabular}%
\end{array}
\label{CBu3}
\end{equation}

The Crank-Nicolson\ scheme is used to discretize time variables of the
unknown $U$ and $V$ in the Coupled Burger equation system which is given (%
\ref{CB1}), we obtain the time discretized form of the equation as%
\begin{equation}
\begin{array}{r}
\dfrac{U^{n+1}-U^{n}}{\Delta t}-\dfrac{U_{xx}^{n+1}+U_{xx}{}^{n}}{2}+k_{1}%
\dfrac{(UU_{x})^{n+1}+(UU_{x}^{n})^{n}}{2}+k_{2}\dfrac{%
(UV)_{x}^{n+1}+(UV)_{x}^{n}}{2}=0 \\ 
\dfrac{V^{n+1}-V^{n}}{\Delta t}-\dfrac{V_{xx}^{n+1}+V_{xx}{}^{n}}{2}+k_{1}%
\dfrac{(VV_{x})^{n+1}+(VV_{x}^{n})^{n}}{2}+k_{3}\dfrac{%
(UV)_{x}^{n+1}+(UV)_{x}^{n}}{2}=0%
\end{array}
\label{CB7}
\end{equation}%
where $U^{n+1}=U(x,t_{n}+\Delta t)$ and $V^{n+1}=V(x,t_{n}+\Delta t).$ The
nonlinear term $(UU_{x})^{n+1},$ $(VV_{x})^{n+1}$ and $(UV)_{x}^{n+1}$in Eq.
(\ref{CB7}) is linearized by using the following form \cite{rubin}:

\begin{equation}
\begin{array}{ll}
(UU_{x})^{n+1} & =U^{n+1}U_{x}^{n}+U^{n}U_{x}^{n+1}-U^{n}U_{x}^{n} \\ 
(VV_{x})^{n+1} & =V^{n+1}V_{x}^{n}+V^{n}V_{x}^{n+1}-V^{n}V_{x}^{n} \\ 
&  \\ 
(UV)_{x}^{n+1} & =(U_{x}V)^{n+1}+(UV_{x})^{n+1} \\ 
& 
=U_{x}^{n+1}V^{n}+U_{x}^{n}V^{n+1}-U_{x}^{n}V^{n}+U^{n+1}V_{x}^{n}+U^{n}V_{x}^{n+1}-U^{n}V_{x}^{n}%
\end{array}
\label{Rubin}
\end{equation}

Substitution the approximate solution (\ref{CBu3}) into (\ref{CB7}) and
evaluating the resulting \ equations at the knots yields the system of the
fully-discretized equations

\begin{equation}
\nu _{m1}\delta _{m-1}^{n+1}+\nu _{m2}\phi _{m-1}^{n+1}+\nu _{m3}\delta
_{m}^{n+1}+\nu _{m4}\phi _{m}^{n+1}+\nu _{m5}\delta _{m+1}^{n+1}+\nu
_{m6}\phi _{m+1}^{n+1}=\nu _{m7}\delta _{m-1}^{n}+\nu _{m8}\delta
_{m}^{n}+\nu _{m9}\delta _{m+1}^{n}  \label{sis1}
\end{equation}

and%
\begin{equation}
\nu _{m10}\delta _{m-1}^{n+1}+\nu _{m11}\phi _{m-1}^{n+1}+\nu _{m12}\delta
_{m}^{n+1}+\nu _{m13}\phi _{m}^{n+1}+\nu _{m14}\delta _{m+1}^{n+1}+\nu
_{m15}\phi _{m+1}^{n+1}=\nu _{m7}\phi _{m-1}^{n}+\nu _{m8}\phi _{m}^{n}+\nu
_{m9}\phi _{m+1}^{n}  \label{sis2}
\end{equation}

where%
\begin{equation*}
\begin{array}{l}
\nu _{m1}=\left( \dfrac{2}{\Delta t}+k_{1}K_{2}+k_{2}L_{2}\right) \alpha
_{1}+\left( k_{1}K_{1}+k_{2}L_{1}\right) \beta _{1}-\gamma _{1} \\ 
\nu _{m2}=\left( k_{1}K_{2}\right) \alpha _{1}+\left( k_{2}K_{1}\right)
\beta _{1} \\ 
\nu _{m3}=\left( \dfrac{2}{\Delta t}+k_{1}K_{2}+k_{2}L_{2}\right) \alpha
_{2}-\gamma _{2} \\ 
\nu _{m4}=\left( k_{1}K_{2}\right) \alpha _{2} \\ 
\nu _{m5}=\left( \dfrac{2}{\Delta t}+k_{1}K_{2}+k_{2}L_{2}\right) \alpha
_{1}-\left( k_{1}K_{1}+k_{2}L_{1}\right) \beta _{1}-\gamma _{1} \\ 
\nu _{m6}=\left( k_{1}K_{2}\right) \alpha _{1}-\left( k_{2}K_{1}\right)
\beta _{1} \\ 
\\ 
\nu _{m7}=\frac{2}{\Delta t}\alpha _{1}+\gamma _{1} \\ 
\nu _{m8}=\frac{2}{\Delta t}\alpha _{2}+\gamma _{2} \\ 
\nu _{m9}=\frac{2}{\Delta t}\alpha _{1}+\gamma _{1} \\ 
\\ 
\nu _{m10}=\left( k_{3}L_{2}\right) \alpha _{1}+\left( k_{3}L_{1}\right)
\beta _{1} \\ 
\nu _{m11}=\left( \frac{2}{\Delta t}+k_{1}L_{2}+k_{3}K_{2}\right) \alpha
_{1}+\left( k_{1}L_{1}+k_{3}K_{1}\right) \beta _{1}-\gamma _{1} \\ 
\nu _{m12}=\left( k_{3}L_{2}\right) \alpha _{2} \\ 
\nu _{m12}=\left( \frac{2}{\Delta t}+k_{1}L_{2}+k_{3}K_{2}\right) \alpha
_{2}-\gamma _{2} \\ 
\nu _{m14}=\left( k_{3}L_{2}\right) \alpha _{1}-\left( k_{3}L_{1}\right)
\beta _{1} \\ 
\nu _{m15}=\left( \frac{2}{\Delta t}+k_{1}L_{2}+k_{3}K_{2}\right) \alpha
_{1}-\left( k_{1}L_{1}+k_{3}K_{1}\right) \beta _{1}-\gamma _{1}%
\end{array}%
\end{equation*}%
\begin{eqnarray*}
&&%
\begin{array}{cc}
K_{1}=\alpha _{1}\delta _{m-1}^{n}+\alpha _{2}\delta _{m}^{n}+\alpha
_{3}\delta _{m+1}^{n} & L_{1}=\alpha _{1}\phi _{m-1}^{n}+\alpha _{2}\phi
_{m}^{n}+\alpha _{3}\phi _{m+1}^{n}%
\end{array}
\\
&&%
\begin{array}{cc}
K_{2}=\beta _{1}\delta _{m-1}^{n}+\beta _{2}\delta _{m}^{n}+\beta _{3}\delta
_{m+1}^{n} & L_{2}=\beta _{1}\phi _{m-1}^{n}+\beta _{2}\phi _{m}^{n}+\beta
_{3}\phi _{m+1}^{n}%
\end{array}%
\end{eqnarray*}%
\begin{eqnarray*}
\alpha _{1} &=&\dfrac{s-ph}{2(phc-s)},\text{ }\alpha _{2}=1, \\
\beta _{1} &=&\dfrac{p(1-c)}{2(phc-s)},\text{ }\beta _{2}=\dfrac{p(c-1)}{%
2(phc-s)} \\
\gamma _{1} &=&\dfrac{p^{2}s}{2(phc-s)},\text{ }\gamma _{2}=-\dfrac{p^{2}s}{%
phc-s}.
\end{eqnarray*}%
The system with (\ref{sis1}) and (\ref{sis2}) can be converted the following
matrices system;%
\begin{equation}
\mathbf{Ax}^{n+1}=\mathbf{Bx}^{n}  \label{CB10}
\end{equation}%
where%
\begin{equation*}
\mathbf{A=}%
\begin{bmatrix}
\nu _{m1} & \nu _{m2} & \nu _{m3} & \nu _{m4} & \nu _{m5} & \nu _{m6} &  & 
&  &  \\ 
\nu _{m10} & \nu _{m11} & \nu _{m12} & \nu _{m13} & \nu _{m14} & \nu _{m15}
&  &  &  &  \\ 
&  & \nu _{m1} & \nu _{m2} & \nu _{m3} & \nu _{m4} & \nu _{m5} & \nu _{m6} & 
&  \\ 
&  & \nu _{m10} & \nu _{m11} & \nu _{m12} & \nu _{m13} & \nu _{m14} & \nu
_{m15} &  &  \\ 
&  &  & \ddots  & \ddots  & \ddots  & \ddots  & \ddots  & \ddots  &  \\ 
&  &  &  & \nu _{m1} & \nu _{m2} & \nu _{m3} & \nu _{m4} & \nu _{m5} & \nu
_{m6} \\ 
&  &  &  & \nu _{m10} & \nu _{m11} & \nu _{m12} & \nu _{m13} & \nu _{m14} & 
\nu _{m15}%
\end{bmatrix}%
\end{equation*}%
\begin{equation*}
\mathbf{B=}%
\begin{bmatrix}
\nu _{m7} & 0 & \nu _{m8} & 0 & \nu _{m9} & 0 &  &  &  &  \\ 
0 & \nu _{m7} & 0 & \nu _{m8} & 0 & \nu _{m9} &  &  &  &  \\ 
&  & \nu _{m7} & 0 & \nu _{m8} & 0 & \nu _{m9} & 0 &  &  \\ 
&  & 0 & \nu _{m7} & 0 & \nu _{m8} & 0 & \nu _{m9} &  &  \\ 
&  &  & \ddots  & \ddots  & \ddots  & \ddots  & \ddots  & \ddots  &  \\ 
&  &  &  & \nu _{m7} & 0 & \nu _{m8} & 0 & \nu _{m9} & 0 \\ 
&  &  &  & 0 & \nu _{m7} & 0 & \nu _{m8} & 0 & \nu _{m9}%
\end{bmatrix}%
\end{equation*}

The system (\ref{CB10}) consist of $2N+2$ linear equation in $2N+6$ unknown
parameters $\mathbf{x}^{n+1}=(\delta _{-1}^{n+1},\phi _{-1}^{n+1},\delta
_{0}^{n+1},\phi _{0}^{n+1},\ldots ,\delta _{n+1}^{n+1},\phi _{n+1}^{n+1},)$.
To obtain a unique solution an additional four constraints\ are needed. By
imposing the Dirichlet boundary conditions this will lead us to the
following relations;%
\begin{equation}
\begin{array}{r}
\delta _{-1}=(f_{1}(a,t)-\alpha _{2}\delta _{0}-\alpha _{3}\delta
_{1})/\alpha _{1} \\ 
\phi _{-1}=(g_{1}(a,t)-\alpha _{2}\phi _{0}-\alpha _{3}\phi _{1})/\alpha _{1}
\\ 
\delta _{N+1}=(f_{2}(b,t)-\alpha _{1}\delta _{N-1}-\alpha _{2}\delta
_{N})/\alpha _{3} \\ 
\phi _{N+1}=(g_{2}(b,t)-\alpha _{1}\phi _{N-1}+\alpha _{2}\phi _{N})/\alpha
_{3}%
\end{array}
\label{CB11}
\end{equation}

\section{The Initial State}

\qquad Initial parameters $\delta _{-1}^{0},$ $\phi _{-1}^{0},$ $\delta
_{0}^{0},$ $\phi _{0}^{0},\ldots ,\delta _{N+1}^{0},$ $\phi _{N+1}^{0}$ can
be determined from the initial condition and first space derivative of the
initial conditions at the boundaries as the following:%
\begin{equation}
\begin{array}{l}
U^{0}(a,0)=\dfrac{s-ph}{2(phc-s)}\delta _{-1}^{0}+\delta _{0}^{0}+\dfrac{s-ph%
}{2(phc-s)}\delta _{1}^{0} \\ 
U^{0}(x_{m},0)=\dfrac{s-ph}{2(phc-s)}\delta _{m-1}^{0}+\delta _{m}^{0}+%
\dfrac{s-ph}{2(phc-s)}\delta _{m+1}^{0},\text{ \  \ }m=1,2,...,N-1 \\ 
U^{0}(b,0)=U_{N}^{0}=\dfrac{s-ph}{2(phc-s)}\delta _{N-1}^{0}+\delta _{N}^{0}+%
\dfrac{s-ph}{2(phc-s)}\delta _{N+1}^{0}%
\end{array}
\label{BU1}
\end{equation}%
and%
\begin{equation}
\begin{array}{l}
V^{0}(a,0)=\dfrac{s-ph}{2(phc-s)}\phi _{-1}^{0}+\phi _{0}^{0}+\dfrac{s-ph}{%
2(phc-s)}\phi _{1}^{0} \\ 
V^{0}(x_{m},0)=\dfrac{s-ph}{2(phc-s)}\phi _{m-1}^{0}+\phi _{m}^{0}+\dfrac{%
s-ph}{2(phc-s)}\phi _{m+1}^{0},\text{ \  \ }m=1,2,...,N-1 \\ 
V^{0}(b,0)=\dfrac{s-ph}{2(phc-s)}\phi _{N-1}^{0}+\phi _{N}^{0}+\dfrac{s-ph}{%
2(phc-s)}\phi _{N+1}^{0}%
\end{array}
\label{BV1}
\end{equation}

The system (\ref{BU1}) which is constituted for initial conditions consists $%
N+1$ equations and $N+3$ unknown, so we have to eliminate $\delta _{-1}^{0}$
and $\delta _{N+1}^{0}$ for solving this system using following derivatives
conditions%
\begin{equation*}
\delta _{-1}^{0}=\delta _{1}^{0}+\frac{2(phc-s)}{p(1-c)}U_{0}^{\prime },%
\text{ }\delta _{N+1}^{0}=\delta _{N-1}^{0}-\frac{2(phc-s)}{p(1-c)}%
U_{N}^{\prime },
\end{equation*}%
and if the equations system is rearranged for the above conditions, then
following form is obtained%
\begin{equation*}
\left[ 
\begin{tabular}{ccccccc}
$\smallskip 1$ & $\tfrac{s-ph}{phc-s}$ &  &  &  &  &  \\ 
$\tfrac{s-ph}{2(phc-s)}$ & $1$ & $\tfrac{s-ph}{2(phc-s)}$ &  &  &  &  \\ 
&  &  & $\ddots $ &  &  &  \\ 
&  &  &  & $\tfrac{s-ph}{2(phc-s)}\smallskip $ & $1$ & $\tfrac{s-ph}{2(phc-s)%
}$ \\ 
&  &  &  &  & $\tfrac{s-ph}{phc-s}$ & $1$%
\end{tabular}%
\  \right] \left[ 
\begin{tabular}{c}
$\delta _{0}^{0}\smallskip $ \\ 
$\delta _{1}^{0}$ \\ 
$\vdots $ \\ 
$\delta _{N-1}^{0}\smallskip $ \\ 
$\delta _{N}^{0}$%
\end{tabular}%
\  \right] =\left[ 
\begin{tabular}{c}
$U_{0}^{\prime }-\tfrac{s-ph}{p(1-c)}U_{0}^{\prime }\smallskip $ \\ 
$U_{1}^{\prime }$ \\ 
$\vdots $ \\ 
$U_{N-1}^{\prime }\smallskip $ \\ 
$U_{N}^{\prime }-\tfrac{s-ph}{p(c-1)}U_{N}^{\prime }$%
\end{tabular}%
\right]
\end{equation*}%
which can also be solved using a variant of the Thomas algorithm. As the
same way, from the system (\ref{BV1}), $\phi _{-1}^{0}$ and $\phi _{N+1}^{0}$
can be eliminated using 
\begin{equation*}
\phi _{-1}^{0}=\phi _{1}^{0}+\frac{2(phc-s)}{p(1-c)}V_{0}^{\prime },\text{ }%
\phi _{N+1}^{0}=\phi _{N-1}^{0}-\frac{2(phc-s)}{p(1-c)}V_{N}^{\prime },
\end{equation*}%
conditions and the following three bounded matrix is obtained.%
\begin{equation*}
\left[ 
\begin{tabular}{ccccccc}
$\smallskip 1$ & $\tfrac{s-ph}{phc-s}$ &  &  &  &  &  \\ 
$\tfrac{s-ph}{2(phc-s)}$ & $1$ & $\tfrac{s-ph}{2(phc-s)}$ &  &  &  &  \\ 
&  &  & $\ddots $ &  &  &  \\ 
&  &  &  & $\tfrac{s-ph}{2(phc-s)}\smallskip $ & $1$ & $\tfrac{s-ph}{2(phc-s)%
}$ \\ 
&  &  &  &  & $\tfrac{s-ph}{phc-s}$ & $1$%
\end{tabular}%
\  \right] \left[ 
\begin{tabular}{c}
$\phi _{0}^{0}\smallskip $ \\ 
$\phi _{1}^{0}$ \\ 
$\vdots $ \\ 
$\phi _{N-1}^{0}\smallskip $ \\ 
$\phi _{N}^{0}$%
\end{tabular}%
\  \right] =\left[ 
\begin{tabular}{c}
$V_{0}^{\prime }-\tfrac{s-ph}{p(1-c)}V_{0}^{\prime }\smallskip $ \\ 
$V_{1}^{\prime }$ \\ 
$\vdots $ \\ 
$V_{N-1}^{\prime }\smallskip $ \\ 
$V_{N}^{\prime }-\tfrac{s-ph}{p(c-1)}V_{N}^{\prime }$%
\end{tabular}%
\right]
\end{equation*}

\section{Numerical Tests}

\qquad In this section, we will compare the efficiency and accuracy of
suggested method problem. The obtained results will compare with \cite%
{2011iki}, \cite{2012}, \cite{2014uc} and \cite{Aksoy}, while $p$ changes.
The accuracy of the schemes is measured in terms of the following discrete
error norm $L_{\infty }$%
\begin{equation*}
L_{\infty }=\left \vert U-U_{N}\right \vert _{\infty }=\max \limits_{j}\left
\vert U_{j}-(U_{N})_{j}^{n}\right \vert .
\end{equation*}%
\textbf{Problem 1) }Consider the Coupled Burgers equation system (\ref{CB1})
with the following initial and boundary conditions%
\begin{equation*}
U(x,0)=\sin (x),\text{ }V(x,0)=\sin (x)
\end{equation*}%
and%
\begin{equation*}
U(-\pi ,t)=U(\pi ,t)=V(-\pi ,t)=V(\pi ,t)=0
\end{equation*}%
The exact solution is%
\begin{equation*}
U(x,t)=V(x,t)=e^{-t}\sin (x)
\end{equation*}

We compute the numerical solutions using the selected values $k_{1}=-2,$ $%
k_{2}=1$ and $k_{3}=1$ with different values of time step length $\Delta t.$%
In our first computation, we take $t=0.1,$ $\Delta t=0.001$ while the number
of partition $N$ changes$.$ The corresponding results are presented in Table
2 a. In our computation, we compute the maximum absolute errors at time
level $t=1$ for the parameters with different decreasing values of $t$. The
corresponding results are reported in Table 2 b. In both computations, the
results are same for $U(x,t)$ and $V(x,t)$ because of symmetric initial and
boundary conditions. And also we correspond the obtained numerical solutions
by different settings of parameters, specifically for those taken by \cite%
{2014uc} in Table 2 c for $N=50$, $\Delta t=0.01$ and increasing $t.$ And
also in Table 2, we present the rate of convergence in space which is
clearly of second order.

\begin{equation*}
\begin{tabular}{|l|}
\hline
Table 2 a: $L_{\infty }$ Error norms for $t=0.1,$ $\Delta t=0.001,$ $%
U(x,t)=V(x,t)$ \\ \hline
\begin{tabular}{ccccc}
& Present ($p=1$) & Present (Various $p$) & \cite{Aksoy}, ($\lambda =0)$ & 
\cite{Aksoy} (Various $\lambda $) \\ 
$N=200$ & $0.01489\times 10^{-5}$ & $0.00121\times 10^{-5}$ & $0.74326\times
10^{-5}$ & $0.00079\times 10^{-5}$ \\ 
&  & $(p=0.00004330000)$ &  & $(\lambda =-1.640\times 10^{-4})$ \\ 
$N=400$ & $0.00372\times 10^{-5}$ & $0.00044\times 10^{-5}$ & $0.18534\times
10^{-5}$ & $0.00006\times 10^{-5}$ \\ 
&  & $(p=0.00012611302)$ &  & $(\lambda =-4.087\times 10^{-5})$%
\end{tabular}
\\ \hline
\end{tabular}%
\end{equation*}%
\begin{equation*}
\begin{tabular}{|l|}
\hline
Table 2 b: $L_{\infty }$ Error norms for $t=1,$ $N=400,$ $U(x,t)=V(x,t)$ \\ 
\hline
\begin{tabular}{ccccc}
& Present ($p=1$) & Present (Various $p$) & \cite{Aksoy}, ($\lambda =0)$ & 
\cite{Aksoy} (Various $\lambda $) \\ 
$\Delta t=0.01$ & $1.8194\times 10^{-5}$ & $0.00247\times 10^{-5}$ & $%
1.08691\times 10^{-5}$ & $0.00131\times 10^{-5}$ \\ 
&  & $(p=0.00010099997)$ &  & $(\lambda =-5.896\times 10^{-5})$ \\ 
$\Delta t=0.001$ & $1.5159\times 10^{-5}$ & $0.00309\times 10^{-5}$ & $%
1.10393\times 10^{-5}$ & $0.00036\times 10^{-5}$ \\ 
&  & $(p=0.00021660000)$ &  & $(\lambda =-5.992\times 10^{-5})$%
\end{tabular}
\\ \hline
\end{tabular}%
\end{equation*}%
\begin{equation*}
\begin{tabular}{|l|}
\hline
Table 2 c: $L_{\infty }$ Error norms for $\Delta t=0.01,$ $N=50,$ $%
U(x,t)=V(x,t)$ different $t.$ \\ \hline
\begin{tabular}{cccccc}
& Present ($p=1$) & Present (Various $p$) & \cite{2012} & \cite{2012iki} & 
\cite{2014uc} \\ 
$t=0.5$ & $7.9881\times 10^{-4}$ & $3.4770\times 10^{-4}$ & $1.51688\times
10^{-4}$ & $2.26627\times 10^{-5}$ & $1.103080984\times 10^{-4}$ \\ 
$t=1.0$ & $9.6837\times 10^{-4}$ & $4.2166\times 10^{-4}$ & $1.83970\times
10^{-4}$ & $1.46179\times 10^{-5}$ & $1.336880384\times 10^{-4}$ \\ 
$t=2.0$ & $7.1154\times 10^{-4}$ & $3.1006\times 10^{-4}$ & $1.35250\times
10^{-4}$ & $0.73805\times 10^{-5}$ & $9.818252567\times 10^{-5}$ \\ 
$t=3.0$ & $3.9213\times 10^{-4}$ & $1.7100\times 10^{-4}$ & $7.46014\times
10^{-4}$ & $0.40272\times 10^{-5}$ & $1.029870405\times 10^{-5}$%
\end{tabular}
\\ \hline
\end{tabular}%
\end{equation*}%
\begin{equation*}
\begin{tabular}{|l|}
\hline
Table 2 d: The rate of convergence for $\Delta t=0.01$ and $\Delta t=0.0001$
respectively \\ \hline
\begin{tabular}{ccc}
$\Delta t=0.01$ & Present ($p=1$) & order \\ 
$N=50$ & $3.9213\times 10^{-4}$ &  \\ 
$N=100$ & $9.9430\times 10^{-5}$ & $1.9773$ \\ 
$N=150$ & $4.4895\times 10^{-5}$ & $1.9649$ \\ 
$N=200$ & $2.5808\times 10^{-5}$ & $1.9192$ \\ 
$N=250$ & $1.6965\times 10^{-5}$ & $1.8734$%
\end{tabular}%
\begin{tabular}{ccc}
$\Delta t=0.0001$ & Present ($p=1$) & order \\ 
$N=50$ & $3.9213\times 10^{-4}$ &  \\ 
$N=100$ & $9.9430\times 10^{-5}$ & $1.9909$ \\ 
$N=150$ & $4.4895\times 10^{-5}$ & $2.0032$ \\ 
$N=200$ & $2.5808\times 10^{-5}$ & $1.9931$ \\ 
$N=250$ & $1.6965\times 10^{-5}$ & $1.9929$%
\end{tabular}
\\ \hline
\end{tabular}%
\end{equation*}%
The corresponding graphical illustrations are presented in Figures 2 for $%
k_{1}=-2,$ $k_{2}=1$, $k_{3}=1,$ $N=400$ and $\Delta t=0.001$ at different $%
t $ for best parameter $p=0.0002166$. In Figure 3-4, computed solutions of $%
v $ different time levels for $k_{1},k_{2}$ fixed and $k_{1},k_{3}$ fixed
respectively.%
\begin{equation*}
\begin{tabular}{l}
\FRAME{itbpF}{2.6818in}{2.2217in}{0in}{}{}{fig2.bmp}{\special{language
"Scientific Word";type "GRAPHIC";maintain-aspect-ratio TRUE;display
"USEDEF";valid_file "F";width 2.6818in;height 2.2217in;depth
0in;original-width 2.6532in;original-height 2.1932in;cropleft "0";croptop
"1";cropright "1";cropbottom "0";filename 'Fig2.bmp';file-properties
"XNPEU";}} \\ 
Figure 2: Numerical Solutions at various \\ 
$t$ for $N=400,$ $\Delta t=0.001,$ $p=0.0002166$%
\end{tabular}%
\end{equation*}%
\begin{equation*}
\begin{array}{c}
\begin{array}{ccc}
\FRAME{itbpF}{2.1223in}{1.7279in}{0in}{}{}{fig3a.bmp}{\special{language
"Scientific Word";type "GRAPHIC";maintain-aspect-ratio TRUE;display
"USEDEF";valid_file "F";width 2.1223in;height 1.7279in;depth
0in;original-width 2.6195in;original-height 2.1266in;cropleft "0";croptop
"1";cropright "1";cropbottom "0";filename 'Fig3a.bmp';file-properties
"XNPEU";}} & \FRAME{itbpF}{2.1664in}{1.7884in}{0in}{}{}{fig3b.bmp}{\special%
{language "Scientific Word";type "GRAPHIC";maintain-aspect-ratio
TRUE;display "USEDEF";valid_file "F";width 2.1664in;height 1.7884in;depth
0in;original-width 2.674in;original-height 2.2001in;cropleft "0";croptop
"1";cropright "1";cropbottom "0";filename 'Fig3b.bmp';file-properties
"XNPEU";}} & \FRAME{itbpF}{2.1664in}{1.7884in}{0in}{}{}{fig3c.bmp}{\special%
{language "Scientific Word";type "GRAPHIC";maintain-aspect-ratio
TRUE;display "USEDEF";valid_file "F";width 2.1664in;height 1.7884in;depth
0in;original-width 2.674in;original-height 2.2001in;cropleft "0";croptop
"1";cropright "1";cropbottom "0";filename 'Fig3c.bmp';file-properties
"XNPEU";}} \\ 
\text{a: }k_{1}=-2,\text{ }k_{2}=1,\text{ }k_{3}=-8 & \text{b: }k_{1}=-2,%
\text{ }k_{2}=1,\text{ }k_{3}=-4 & \text{c: }k_{1}=-2,\text{ }k_{2}=1,\text{ 
}k_{3}=-0%
\end{array}
\\ 
\text{Figure 3: Computed solutions of }V\text{ Problem 1 for different time
levels (}k_{1},k_{2}\text{ fixed)}%
\end{array}%
\end{equation*}%
\begin{equation*}
\begin{array}{c}
\begin{array}{ccc}
\FRAME{itbpF}{2.1612in}{1.7884in}{0in}{}{}{fig4a.bmp}{\special{language
"Scientific Word";type "GRAPHIC";maintain-aspect-ratio TRUE;display
"USEDEF";valid_file "F";width 2.1612in;height 1.7884in;depth
0in;original-width 2.6671in;original-height 2.2001in;cropleft "0";croptop
"1";cropright "1";cropbottom "0";filename 'Fig4a.bmp';file-properties
"XNPEU";}} & \FRAME{itbpF}{2.1767in}{1.74in}{0in}{}{}{fig4b.bmp}{\special%
{language "Scientific Word";type "GRAPHIC";maintain-aspect-ratio
TRUE;display "USEDEF";valid_file "F";width 2.1767in;height 1.74in;depth
0in;original-width 2.687in;original-height 2.1395in;cropleft "0";croptop
"1";cropright "1";cropbottom "0";filename 'Fig4b.bmp';file-properties
"XNPEU";}} & \FRAME{itbpF}{2.1664in}{1.7884in}{0in}{}{}{fig4c.bmp}{\special%
{language "Scientific Word";type "GRAPHIC";maintain-aspect-ratio
TRUE;display "USEDEF";valid_file "F";width 2.1664in;height 1.7884in;depth
0in;original-width 2.674in;original-height 2.2001in;cropleft "0";croptop
"1";cropright "1";cropbottom "0";filename 'Fig4c.bmp';file-properties
"XNPEU";}} \\ 
\text{a: }k_{1}=-2,\text{ }k_{2}=-8,\text{ }k_{3}=1 & \text{b: }k_{1}=-2,%
\text{ }k_{2}=-4,\text{ }k_{3}=1 & \text{c: }k_{1}=-2,\text{ }k_{2}=0,\text{ 
}k_{3}=1%
\end{array}
\\ 
\text{Figure 4: Computed solutions of }V\text{ Problem 1 for different time
levels (}k_{1},k_{3}\text{ fixed)}%
\end{array}%
\end{equation*}%
\begin{equation*}
\begin{array}{c}
\begin{array}{ccc}
\FRAME{itbpF}{2.1612in}{1.7997in}{0in}{}{}{fig5a.bmp}{\special{language
"Scientific Word";type "GRAPHIC";maintain-aspect-ratio TRUE;display
"USEDEF";valid_file "F";width 2.1612in;height 1.7997in;depth
0in;original-width 2.6671in;original-height 2.2139in;cropleft "0";croptop
"1";cropright "1";cropbottom "0";filename 'Fig5a.bmp';file-properties
"XNPEU";}} & \FRAME{itbpF}{2.1664in}{1.7997in}{0in}{}{}{fig5b.bmp}{\special%
{language "Scientific Word";type "GRAPHIC";maintain-aspect-ratio
TRUE;display "USEDEF";valid_file "F";width 2.1664in;height 1.7997in;depth
0in;original-width 2.674in;original-height 2.2139in;cropleft "0";croptop
"1";cropright "1";cropbottom "0";filename 'Fig5b.bmp';file-properties
"XNPEU";}} & \FRAME{itbpF}{2.1612in}{1.7461in}{0in}{}{}{fig5c.bmp}{\special%
{language "Scientific Word";type "GRAPHIC";maintain-aspect-ratio
TRUE;display "USEDEF";valid_file "F";width 2.1612in;height 1.7461in;depth
0in;original-width 2.6671in;original-height 2.1465in;cropleft "0";croptop
"1";cropright "1";cropbottom "0";filename 'Fig5c.bmp';file-properties
"XNPEU";}} \\ 
\text{a: }k_{1}=-8,\text{ }k_{2}=1,\text{ }k_{3}=1 & \text{b: }k_{1}=-4,%
\text{ }k_{2}=1,\text{ }k_{3}=1 & \text{c: }k_{1}=0,\text{ }k_{2}=1,\text{ }%
k_{3}=1%
\end{array}
\\ 
\text{Figure 5: Computed solutions of }V\text{ Problem 1 for different time
levels (}k_{2},k_{3}\text{ fixed)}%
\end{array}%
\end{equation*}%
\textbf{\ Problem 2) }Numerical solutions of considered coupled Burgers'
equations are obtained for $k_{1}=2$ with different values of $k_{2}$ and $%
k_{3}$ at different time levels. In this situation the exact solution is%
\begin{equation*}
\begin{array}{l}
U(x,t)=a_{0}-2A(\frac{2k_{2}-1}{4k_{2}k_{3}-1})\tanh (A(x-2At)) \\ 
\\ 
V(x,t)=a_{0}(\frac{2k_{3}-1}{2k_{2}-1})-2A(\frac{2k_{2}-1}{4k_{2}k_{3}-1}%
)\tanh (A(x-2At))%
\end{array}%
\end{equation*}

Thus, the initial and boundary conditions are taken from the exact solution
is%
\begin{equation*}
\begin{array}{l}
U(x,0)=a_{0}-2A(\frac{2k_{2}-1}{4k_{2}k_{3}-1})\tanh (Ax) \\ 
\\ 
V(x,0)=a_{0}(\frac{2k_{3}-1}{2k_{2}-1})-2A(\frac{2k_{2}-1}{4k_{2}k_{3}-1}%
)\tanh (Ax)%
\end{array}%
\end{equation*}

Thus, the initial and boundary conditions are exracted from the exact
solution. Where $a_{0}=0.05$ and $A=\dfrac{1}{2}(\frac{a_{0}(4k_{2}k_{3}-1)}{%
2k_{2}-1}).$The numerical solutions have been computed for the domain $x\in
\lbrack 0,1],$ $\Delta t=0.001$ and number of partitions $N=10$ and $100.$
The maximum absolute errors have been computed and compared in Tables 3 a-3
b for $t=1$ with those available in the literature \cite{Aksoy}.%
\begin{equation*}
\begin{tabular}{|l|}
\hline
Table 3 a: $L_{\infty }$ Error norms for $t=1,$ $\Delta t=0.001,$ $U(x,t),$ $%
k_{1}=2,$ $k_{2}=1$ and $k_{3}=0.3$ \\ \hline
\begin{tabular}{cccc}
& Present ($p=1$) & \cite{Aksoy}, ($\lambda =0)$ & \cite{Aksoy} (Various $%
\lambda $) \\ 
$N=10$ & $3.7323\times 10^{-6}$ & $3.73505\times 10^{-5}$ & 
\multicolumn{1}{l}{$0.00077\times 10^{-5}(\lambda =6\times 10^{-5})$} \\ 
$N=100$ & $3.7350\times 10^{-6}$ & $3.73503\times 10^{-5}$ & 
\multicolumn{1}{l}{$0.00078\times 10^{-5}(\lambda =-4.087\times 10^{-5})$}%
\end{tabular}
\\ \hline
\end{tabular}%
\end{equation*}%
\begin{equation*}
\begin{tabular}{|l|}
\hline
Table 3 b: $L_{\infty }$ Error norms for $t=1,$ $\Delta t=0.001,$ $V(x,t),$ $%
k_{1}=2,$ $k_{2}=1$ and $k_{3}=0.3$ \\ \hline
\begin{tabular}{cccc}
& Present ($p=1$) & \cite{Aksoy}, ($\lambda =0)$ & \cite{Aksoy} (Various $%
\lambda $) \\ 
$N=10$ & $1.2569\times 10^{-6}$ & $1.29030\times 10^{-5}$ & 
\multicolumn{1}{l}{$0.00079\times 10^{-5}(\lambda =-6\times 10^{-5})$} \\ 
$N=100$ & $1.2871\times 10^{-6}$ & $1.29038\times 10^{-5}$ & 
\multicolumn{1}{l}{$0.00079\times 10^{-5}(\lambda =-4.087\times 10^{-4})$}%
\end{tabular}
\\ \hline
\end{tabular}%
\end{equation*}%
\begin{equation*}
\begin{tabular}{|l|}
\hline
Table 3 c: Maximum error norms for $U(x,t)$ in Problem 2 ($N=21,\Delta
t=0.01,$ $k_{1}=2)$ \\ \hline
\begin{tabular}{cccccc}
$t$ & $k_{2}$ & $k_{3}$ & Present ($p=1$) & \cite{2012} & \cite{2014uc} \\ 
$0.5$ & $0.1$ & $0.3$ & $8.8160\times 10^{-6}$ & $4.173\times 10^{-5}$ & $%
4.189217417\times 10^{-5}$ \\ 
& $0.3$ & $0.03$ & $9.2556\times 10^{-6}$ & $4.585\times 10^{-5}$ & $%
4.584830094\times 10^{-5}$ \\ 
$1.0$ & $0.1$ & $0.3$ & $8.8878\times 10^{-6}$ & $8.275\times 10^{-5}$ & $%
8.269641708\times 10^{-5}$ \\ 
& $0.3$ & $0.03$ & $9.3324\times 10^{-6}$ & $9.167\times 10^{-5}$ & $%
9.147335667\times 10^{-5}$ \\ 
$3.0$ & $0.1$ & $0.3$ & $8.9174\times 10^{-6}$ & $2.408\times 10^{-4}$ & $%
2.401202768\times 10^{-4}$ \\ 
& $0.1$ & $0.03$ & $9.3691\times 10^{-6}$ & $2.747\times 10^{-4}$ & $%
2.704203611\times 10^{-4}$%
\end{tabular}
\\ \hline
\end{tabular}%
\end{equation*}%
\begin{equation*}
\begin{tabular}{|l|}
\hline
Table 3 d: Maximum error norms for $V(x,t)$ in Problem 2 ($N=21,\Delta
t=0.01,$ $k_{1}=2)$ \\ \hline
\begin{tabular}{cccccc}
$t$ & $k_{2}$ & $k_{3}$ & Present ($p=1$) & \cite{2012} & \cite{2014uc} \\ 
$0.5$ & $0.1$ & $0.3$ & $2.8380\times 10^{-6}$ & $5.418\times 10^{-5}$ & $%
9.094743099\times 10^{-6}$ \\ 
& $0.3$ & $0.03$ & $1.1179\times 10^{-5}$ & $2.826\times 10^{-5}$ & $%
2.48218881\times 10^{-5}$ \\ 
$1.0$ & $0.1$ & $0.3$ & $2.8686\times 10^{-6}$ & $1.074\times 10^{-4}$ & $%
1.696286567\times 10^{-5}$ \\ 
& $0.3$ & $0.03$ & $1.1269\times 10^{-5}$ & $5.673\times 10^{-5}$ & $%
4.965329678\times 10^{-5}$ \\ 
$3.0$ & $0.1$ & $0.3$ & $2.9081\times 10^{-6}$ & $3.119\times 10^{-4}$ & $%
4.505480184\times 10^{-5}$ \\ 
& $0.1$ & $0.03$ & $1.1301\times 10^{-5}$ & $1.663\times 10^{-4}$ & $%
1.498311672\times 10^{-5}$%
\end{tabular}
\\ \hline
\end{tabular}%
\end{equation*}

\begin{equation*}
\begin{tabular}{c}
$%
\begin{array}{cc}
\FRAME{itbpF}{3.2283in}{2.668in}{0in}{}{}{fig6a.bmp}{\special{language
"Scientific Word";type "GRAPHIC";maintain-aspect-ratio TRUE;display
"USEDEF";valid_file "F";width 3.2283in;height 2.668in;depth
0in;original-width 2.6671in;original-height 2.2001in;cropleft "0";croptop
"1";cropright "1";cropbottom "0";filename 'Fig6a.bmp';file-properties
"XNPEU";}} & \FRAME{itbpF}{3.2283in}{2.668in}{0in}{}{}{fig6b.bmp}{\special%
{language "Scientific Word";type "GRAPHIC";maintain-aspect-ratio
TRUE;display "USEDEF";valid_file "F";width 3.2283in;height 2.668in;depth
0in;original-width 2.6671in;original-height 2.2001in;cropleft "0";croptop
"1";cropright "1";cropbottom "0";filename 'Fig6b.bmp';file-properties
"XNPEU";}}%
\end{array}%
$ \\ 
Figure 6: Numerical Solutions for $U(x,t)$ and $V(x,t)$, $N=21,$ $\Delta
t=0.001,$ $t=1,$ $0\leq x\leq 1,$ $k_{2}=0.1$ and $k_{3}=0.3$%
\end{tabular}%
\end{equation*}%
\textbf{Problem 3)\ }Consider the Coupled Burger Equation system (\ref{CB1})
with the following initial conditions

\begin{equation*}
U(x,0)=\left \{ 
\begin{array}{cc}
\sin (2\pi x), & x\in \lbrack 0,0.5] \\ 
0, & x\in (0.5,1]%
\end{array}%
\right.
\end{equation*}

\begin{equation*}
V(x,0)=\left \{ 
\begin{array}{cc}
0, & x\in \lbrack 0,0.5] \\ 
-\sin (2\pi x), & x\in (0.5,1]%
\end{array}%
\right.
\end{equation*}%
and zero boundary conditions. In the Problem 3, the solutions have been
carried out on $x\in \lbrack 0,1]$ with $\Delta t=0.001$ and number of
partitions as $50$. Maximum values of $u$ and $v$ at different time levels
for $k_{2}=k_{3}=10$ have been given in Table 4a and 4 b, while the Tables 4
c and 4 d represent the maximum values for $k_{2}=k_{3}=100.$%
\begin{equation*}
\begin{tabular}{|l|}
\hline
Table 4 a: Maximum values of $U$ at different time levels for $%
k_{2}=k_{3}=10 $ \\ \hline
\begin{tabular}{ccccc}
$t$ & Present ($p=1$) & \cite{2011iki} & \cite{2014uc} & at point \\ 
$0.1$ & $0.144501$ & $0.14456$ & $0.144491495800$ & $0.58$ \\ 
$0.2$ & $0.052353$ & $0.05237$ & $0.052356151890$ & $0.54$ \\ 
$0.3$ & $0.019317$ & $0.01932$ & $0.019318838080$ & $0.52$ \\ 
$0.4$ & $0.007183$ & $0.00718$ & $0.007184856672$ & $0.50$%
\end{tabular}
\\ \hline
\end{tabular}%
\end{equation*}%
\begin{equation*}
\begin{tabular}{|l|}
\hline
Table 4 b: Maximum values of $V$ at different time levels for $%
k_{2}=k_{3}=10 $ \\ \hline
\begin{tabular}{ccccc}
$t$ & Present ($p=1$) & \cite{2011iki} & \cite{2014uc} & at point \\ 
$0.1$ & $0.143155$ & $0.14306$ & $0.143141957500$ & $0.66$ \\ 
$0.2$ & $0.047003$ & $0.04697$ & $0.047006446750$ & $0.56$ \\ 
$0.3$ & $0.017258$ & $0.01725$ & $0.017260356430$ & $0.52$ \\ 
$0.4$ & $0.006415$ & $0.00641$ & $0.006416614856$ & $0.50$%
\end{tabular}
\\ \hline
\end{tabular}%
\end{equation*}%
\begin{equation*}
\begin{tabular}{|l|}
\hline
Table 4 c: Maximum values of $U$ at different time levels for $%
k_{2}=k_{3}=100$ \\ \hline
\begin{tabular}{ccccc}
$t$ & Present ($p=1$) & \cite{2011iki} & \cite{2014uc} & at point \\ 
$0.1$ & $0.04168$ & $0.04175$ & $0.041682987260$ & $0.46$ \\ 
$0.2$ & $0.01476$ & $0.01479$ & $0.014770415340$ & $0.58$ \\ 
$0.3$ & $0.00533$ & $0.00534$ & $0.005337325631$ & $0.54$ \\ 
$0.4$ & $0.00197$ & $0.00198$ & $0.001978065014$ & $0.52$%
\end{tabular}
\\ \hline
\end{tabular}%
\end{equation*}%
\begin{equation*}
\begin{tabular}{|l|}
\hline
Table 4 d: Maximum values of $V$ at different time levels for $%
k_{2}=k_{3}=100$ \\ \hline
\begin{tabular}{ccccc}
$t$ & Present ($p=1$) & \cite{2011iki} & \cite{2014uc} & at point \\ 
$0.1$ & $0.05074$ & $0.05065$ & $0.050737669860$ & $0.76$ \\ 
$0.2$ & $0.01035$ & $0.01033$ & $0.010356602970$ & $0.64$ \\ 
$0.3$ & $0.00351$ & $0.00350$ & $0.003517189432$ & $0.56$ \\ 
$0.4$ & $0.00129$ & $0.00129$ & $0.001294450199$ & $0.52$%
\end{tabular}
\\ \hline
\end{tabular}%
\end{equation*}

Figs. 7, 8 and 9 show the numerical results obtained for different time
levels $t\in \lbrack 0,1]$ at $k_{2}=k_{3}=10$ for $U$ and $V$ with
different values of $k_{1}$. From the Figs. 7-9, it can be easily seen that
the numerical solutions $U^{n}$ and $V^{n}$ decay to zero as $t$ and $k_{1}$
increased.%
\begin{equation*}
\begin{array}{c}
\begin{array}{cc}
\FRAME{itbpF}{2.5633in}{2.3739in}{0in}{}{}{fig7a.bmp}{\special{language
"Scientific Word";type "GRAPHIC";maintain-aspect-ratio TRUE;display
"USEDEF";valid_file "F";width 2.5633in;height 2.3739in;depth
0in;original-width 5.073in;original-height 4.6933in;cropleft "0";croptop
"1";cropright "1";cropbottom "0";filename 'Fig7a.bmp';file-properties
"XNPEU";}} & \FRAME{itbpF}{2.5633in}{2.3739in}{0in}{}{}{fig7b.bmp}{\special%
{language "Scientific Word";type "GRAPHIC";maintain-aspect-ratio
TRUE;display "USEDEF";valid_file "F";width 2.5633in;height 2.3739in;depth
0in;original-width 5.073in;original-height 4.6933in;cropleft "0";croptop
"1";cropright "1";cropbottom "0";filename 'Fig7b.bmp';file-properties
"XNPEU";}}%
\end{array}
\\ 
\text{Fig 7: Num. Sol. }U(x,t)\text{ and }V(x,t)\text{ of Problem 3 at
different time levels for }k_{2}=k_{3}=10\text{ while }k_{1}=1%
\end{array}%
\end{equation*}%
\begin{equation*}
\begin{array}{c}
\begin{array}{cc}
\FRAME{itbpF}{2.5417in}{2.3739in}{0.0104in}{}{}{fig8a.bmp}{\special{language
"Scientific Word";type "GRAPHIC";maintain-aspect-ratio TRUE;display
"USEDEF";valid_file "F";width 2.5417in;height 2.3739in;depth
0.0104in;original-width 5.0272in;original-height 4.6933in;cropleft
"0";croptop "1";cropright "1";cropbottom "0";filename
'Fig8a.bmp';file-properties "XNPEU";}} & \FRAME{itbpF}{2.5417in}{2.3739in}{%
0in}{}{}{fig8b.bmp}{\special{language "Scientific Word";type
"GRAPHIC";maintain-aspect-ratio TRUE;display "USEDEF";valid_file "F";width
2.5417in;height 2.3739in;depth 0in;original-width 5.0272in;original-height
4.6933in;cropleft "0";croptop "1";cropright "1";cropbottom "0";filename
'Fig8b.bmp';file-properties "XNPEU";}}%
\end{array}
\\ 
\text{Fig 8: Num. Sol. }U(x,t)\text{ and }V(x,t)\text{ of Problem 3 at
different time levels for }k_{2}=k_{3}=10\text{ while }k_{1}=10%
\end{array}%
\end{equation*}%
\begin{equation*}
\begin{array}{c}
\begin{array}{cc}
\FRAME{itbpF}{2.5417in}{2.3739in}{0in}{}{}{fig9a.bmp}{\special{language
"Scientific Word";type "GRAPHIC";maintain-aspect-ratio TRUE;display
"USEDEF";valid_file "F";width 2.5417in;height 2.3739in;depth
0in;original-width 5.0272in;original-height 4.6933in;cropleft "0";croptop
"1";cropright "1";cropbottom "0";filename 'Fig9a.bmp';file-properties
"XNPEU";}} & \FRAME{itbpF}{2.5417in}{2.3739in}{0in}{}{}{fig9b.bmp}{\special%
{language "Scientific Word";type "GRAPHIC";maintain-aspect-ratio
TRUE;display "USEDEF";valid_file "F";width 2.5417in;height 2.3739in;depth
0in;original-width 5.0272in;original-height 4.6933in;cropleft "0";croptop
"1";cropright "1";cropbottom "0";filename 'Fig9b.bmp';file-properties
"XNPEU";}}%
\end{array}
\\ 
\text{Fig 9: Num. Sol. }U(x,t)\text{ and }V(x,t)\text{ of Problem 3 at
different time levels for }k_{2}=k_{3}=10\text{ while }k_{1}=50%
\end{array}%
\end{equation*}

\section{Conclusion}

The collocation method together with the exponential B-spline as trial
functions has presented to get the numerical solutions of the coupled
Burgers' equation system. The free parameter in the exponential B-splines is
searched experimentally to get the best numerical solution for the first
problem. In the other problem, results are documented for $p=1$ as an
example. The proposed method has produced less error than the methods listed
in the tables for some text problems. The results are satisfactory and
competent with some available solutions in the literature. Another
advantages is that the method can be used without the complex calculations.
to solve the system of differential equations reliably. And also the
collocation method together with B-spline approximations represents an
economical alternative since it only requires the evaluation of the unknown
parameters at the grid points.

\end{document}